\documentclass[10pt]{article}
\usepackage[cp1251]{inputenc}
\usepackage[english,russian]{babel}     
\usepackage{amsfonts}
\usepackage{amsthm} 
\usepackage{amssymb}
\usepackage{latexsym} 

\newtheorem{teo}{\quad Теорема}

\newtheorem{sled}{\quad Следствие}
\newtheorem{df}{\quad Определение}

\newtheorem{vspteo}{\quad Теорема}
\newtheorem{vsplem}{\quad Лемма}

\renewcommand{\refname}{Reference}

\begin{document}
\Large

\noindent УДК 517.5

\bigskip\noindent
{\bf С.О. Чайченко}
( Донбасский государственный педагогический университет, Славянск, Украина)

\bigskip\noindent
{\bf  ПРИБЛИЖЕНИЯ  ПЕРИОДИЧЕСКИХ ФУНКЦИЙ \\ В ВЕСОВЫХ ПРОСТРАНСТВАХ ОРЛИЧА}

\bigskip
{\noindent В работе получен ряд  прямых и обратных теорем теории приближения для $\psi$-дифференцируемых функций в  метриках весовых пространств Орлича с весами, которые принадлежат классу Макенхаупа. }
\bigskip

{ \bf Ключевые слова:} пространства Орлича, неравенство Бернштейна, наилучшее приближение, $\psi$-производная, прямые и обратные теоремы, классы Макенхаупа.
\bigskip

{ \bf AMS:} 46E30, 42A10,41A17,41A20,41A25,41A27, 41A30.

\bigskip
\bigskip
\bigskip

{\bf 1. Определения и постановка задачи.} Приведем некоторые сведения из теории выпуклых функций и весовых пространств Орлича (см. \cite{Geneb_Gogat_Kokilash_Krbec, Krasnosel_Rutick_Vipuk_funkt}).

\begin{df}
Нерерывная выпуклая функция $\Phi=\Phi(x)$ называется функцией Юнга, если $\Phi$ является четной и удовлетворяет условиям
$$
\lim\limits_{x\to 0} \frac{\Phi(x)}{x}=0, \quad \lim\limits_{x\to \infty} \frac{\Phi(x)}{x}=\infty.
$$
\end{df}

\begin{df}
Говорят, что функция  $\Phi$ удовлетворяет условию $\Delta_2~ (\Phi \in \Delta_2)$, если существует постоянная $ c>0$ такая, что
$$
    \Phi(2x)\le c ~\Phi(x), \quad \forall x \in \mathbb{R}.
$$
\end{df}

\begin{df}
Неотрицательная функция $M=M(t),~ t\geq0$ называется квазивыпуклой функцией Юнга, если существует выпуклая функция Юнга $\Phi$ и постоянная $ c>1$ такая, что выполняется неравенство

$$
\Phi(x)\leq M(x)\leq\Phi(cx), \quad \forall x\geq 0.
$$
\end{df}

Обозначим через ${\cal Q}$ множество всех квазивыпуклых функций Юнга.
\begin{df}
 Пусть $M \in {\cal Q}.$ Тогда через $\tilde L_{M,\omega}$ обозначают класс $2\pi$-периодических измеримых по Лебегу функций $f$, которые удовлетворяют условию
$$
    \int\limits^{2\pi}_0 M(|f(x)|)\omega(x)~dx<\infty,
$$
где $\omega(x)$ --- $2\pi$-периодичная измеримая  и почти везде положительная функция (вес), а через $L_{M,\omega}$ обозначают линейную оболочку класса $\tilde L_{M \omega}$.
\end{df}

Множество $\tilde L_{M \omega}$ становиться нормированным пространством, если
$$
    \|f\|_{M,\omega}:=\sup \Bigg\{\int\limits_0^{2\pi} |f(t)g(t)| \omega(t)~dt: \quad
    \int\limits^{2\pi}_0 \tilde M (|g(t)|)\omega(t)~dt \leq 1 \Bigg\},
$$
где $\tilde M(y):=\sup \limits_{x \geq 0} (xy-M(x)),~ y \geq 0$ --- дополнительная в смысле Юнга функция.

Говорят, что весовая функция  $\omega=\omega(t)$ принадлежит классу Макенхаупа $A_p$, $1<p<\infty$, если $\omega$ является $2\pi$-периодической и
$$
    \Bigg(\frac{1}{b-a}\int\limits_a^b\omega(t)~dt \Bigg)
    \Bigg( \frac{1}{b-a} \int\limits^b_a \frac {1}{ \omega^{1/(p-1)}(t)}~dt\Bigg)^{p-1}
    \leq c={\rm const},
$$
где $[a,b]$ произвольный отрезок из $[0,2\pi]$.

Для квазивыпуклой функции $M$ определим величину
$$
    \frac{1}{p(M)}:=\inf\{p: \quad p>0,~ M^p \in {\cal Q}\},
$$
$$
    p'(M):=\frac{p(M)}{p(M)-1},
$$
которая впервые была введена в работе \cite{Gogat.Kokil.1994}. Если $\omega \in A_{p(M)},$ то
$L_{M,\omega} \in L,$ где $L$ --- пространство $2\pi$-периодических суммируемых на периоде
функций и $L_{M,\omega}$ становиться банаховым пространством с нормой Орлича. Банахово пространство
$L_{M,\omega}$ называется весовым пространством Орлича.

Через ${\cal Q}_2^\theta$ обозначают класс функций $M=M(t)$, которые удовлетворяют условию $\Delta_2$  и таких, что $M^\theta$ является  квазивыпуклой для произвольного $\theta\in(0;1)$.

Пусть $f \in L_{M, \omega},$ $M \in {\cal Q}_2^\theta$, $\omega \in A_{p(M)}$  и
$$
    (\sigma_h f):=\frac{1}{2h} \int_{-h}^h f(x+t)~dt, \quad 0<h<\pi, \quad x \in [-\pi;\pi]
$$
--- оператор сдвига. Обозначив через $ I$ --- единичный оператор, рассмотрим величины
$$
    \Omega_r (f; \delta)_{M, \omega}:=
    \sup_{0<h_i<\delta, 1\le i\le r} \bigg\|\prod_{i=1}^{r}
    (I-\sigma_{h_i}) f \bigg\|_{M,\omega}, \quad (\delta>0,~r=1,2,\ldots,)
$$
которые называются модулями гладкости порядка $r$ функции $f.$

Известно (см., например, \cite{Akgun-Weighted_Smirnov_Orlicz}), что модули гладкости $\Omega_r (f; \delta)_{M, \omega}$ обладают такими свойствами:
\begin{enumerate}
    \item Величина $\Omega_r (f; \delta)_{M, \omega}$ является неотрицательной и не убывает, как функция переменной $\delta>0.$
    \item Справедливо неравенство
$$
    \Omega_r (f_1+f_2; \delta)_{M, \omega} \le \Omega_r (f_1; \delta)_{M, \omega}+
    \Omega_r (f_2; \delta)_{M, \omega}.
$$
    \item Если $f \in L_{M, \omega},$ $M \in {\cal Q}_2^\theta$, $\omega \in A_{p(M)}$, то
$$
    \Omega_r (f; \delta)_{M, \omega} \le C \|f\|_{M,\omega}.
$$
\end{enumerate}

Далее нам понадобятся определения $\psi$-интеграла и $\psi$-производной, которые принадлежат А.И.~Степанцу \cite[c.~149 -- 150]{Stepanetz-2002m1} (см. также \cite[c.~137 -- 138]{Stepanetz-2002-angl}).

\begin{df}
Пусть $f\in L$ и
\begin{equation}\label{1}
    S[f]=\frac{a_0(f)}{2}+\sum_{k=1}^\infty(a_k(f)\cos~kx+b_k(f)\sin~kx)
    \equiv \sum_{k=0}^\infty A_k(f,x)
\end{equation}
--- ряд Фурье функции $f$. Пусть, далее,  $\psi(k)=(\psi_1; \psi_2)$
--- пара произвольных числовых последовательностей  $\psi_1(k)$ и $\psi_2(k),~ k=1,2,\ldots~.$.
Рассмотрим ряд
\begin{equation}\label{2}
    A_0+\sum_{k=1}^{\infty}(\psi_1(k)A_k(f,x)+\psi_2(k)\tilde{A}_k(f,x)),
\end{equation}
где $A_0$ --- некоторое число и
$$
    \tilde{A}_k(f,x)=a_k \sin kx-b_k \cos kx.
$$
Если ряд (\ref{2}) для данной функции $f$ и пары $\psi$ является рядом Фурье некоторой функции из $F\in L,$ то  функцию $F$ называют  $\psi$-интегралом функции $f$ и обозначают $F(\cdot)={\cal J}^\psi(f;\cdot).$
Множество $\psi$-интегралов всех функций из $L$ обозначается $L^\psi.$
\end{df}

\begin{df}
Пусть $f\in L$, (\ref{1}) --- ее ряд Фурье и пара
$\psi=(\psi_1,\psi_2)$ удовлетворяет условию
\begin{equation}\label{2*}
    \psi^2(k)=\psi_1^2(k)+\psi_2^2(k)\not=0,
    \quad k\in \mathbb N.
\end{equation}
Если ряд
\begin{equation}\label{3}
    \sum_{k=1}^{\infty}\Bigg(\frac{\psi_1(k)}{\psi^2(k)}
    A_k(f;x)-\frac{\psi_2(k)}{\psi^2(k)}\tilde{A}_k(f;x)\Bigg)
\end{equation}
является рядом Фурье некоторой функции  $\varphi \in L$,   то
$\varphi$ назовем $\psi$-производной функции  $f$ и будем
писать $\varphi(\cdot)=D^{\psi}(f;\cdot)$
$=f^{\psi}(\cdot)$.

Подмножество функций $f\in L$, у которых существуют
$\psi$-производные, обозначают через $\bar{L}^{\psi}$.
\end{df}

Связь между $\psi$-интегралами и  $\psi$-производными установлена в монографии \cite[c.~150]{Stepanetz-2002m1} (см. также \cite[c.~138 -- 139]{Stepanetz-2002-angl}). Там показано, что если функция $f\in L$, ряд (\ref{1}) --- ее ряд Фурье, и выполнено условие  (\ref{2*}), то функция ${\cal J}^{\psi}(f;x)$
обладает $\psi$-производной и справедливо равенство
$$
    D^{\psi}({\cal J}^{\psi}(f;\cdot))=
    f(\cdot)-\frac{a_0}{2}.
$$
Если же функция $f\in \bar{L}^{\psi}$, и ряд (\ref{1})
--- ее ряд Фурье, то функция $D^{\psi}(f;x)$ обладает
$\psi$-интегралом и при этом
$$
    {\cal J}^{\psi}(D^{\psi}(f;{\cdot}))=f(\cdot)+A_0,
$$
где $A_0$ --- некоторая постоянная.

Отметим также, что в случае, когда
$$
    \psi_1(k)=\psi(k)\cos \frac{\beta \pi}{2}, \quad
    \psi_2(k)=\psi(k)\sin \frac{\beta \pi}{2},
$$
определение $\psi$-производной совпадает с определением $(\psi;\beta)$-производной \cite[c.~132]{Stepanetz-2002m1} (см. также \cite[c.~120]{Stepanetz-2002-angl}). Если к тому же $\psi(k)=k^{-r},~r>0,$ то понятие $\psi$-производной совпадает с хорошо известным понятием $(r;\beta)$-производной в смысле Вейля-Надя, которая, в свою очередь, является обобщением дробной производной Вейля и обычной производной порядка $r.$

\bigskip

{\bf 2. Вспомогательные результаты.} При доказательстве основных результатов этой работы будет использоваться факт ограниченности в пространствах $L_{M,\omega}$ оператора Фурье и оператора тригонометрического сопряжения, который был установлен в \cite[c.~278]{Geneb_Gogat_Kokilash_Krbec}.

\begin{vspteo}\label{T.A}
Если $M \in {\cal Q}_2^\theta$ и $\omega \in A_{p(M)},$  то для произвольной функции
$f \in L_{M,\omega}$ выполняются неравенства
\begin{equation}\label{5}
    \| S_n(f) \|_{M, \omega} \le  C \| f \|_{M, \omega},
\end{equation}
и
\begin{equation}\label{6}
    \| \tilde{f}\|_{M, \omega} \le  C \| f \|_{M, \omega},
\end{equation}
где
\begin{equation}\label{9}
    S_n(f; x)=\frac{a_0(f)}{2}+\sum_{k=1}^{n}
    (a_k(f)\cos~kx+b_k(f)\sin~kx), \quad
    n=0,1,\ldots,
\end{equation}
--- частные суммы Фурье порядка $n$ функции $f$, $\tilde{f}(\cdot)$ --- функция, тригонометрически сопряженная с $f(\cdot)$ и $C>0$ --- некоторая постоянная, которая не зависит от $f$ и $n$.
\end{vspteo}

Обозначим через
$$
    E_{n} (\varphi)_{M, \omega}:=\inf\limits_{t_{n-1} \in
    {\cal T}_{n-1}} \| \varphi- t_{n-1}\|_{M, \omega},
    \quad \varphi \in L_{M,\omega},
$$
--- наилучшее приближение функции $\varphi$ с помощью подпространства ${\cal T}_{n-1}$
тригонометрических полиномов порядка, не выше $n-1.$ Из теоремы \ref{T.A}, в частности, следует, что
\begin{equation}\label{7}
    \| f-S_{n-1}(f) \|_{M, \omega} = {\cal O}(1) E_{n}(f)_{M, \omega} =
    {\cal O}(1)E_{n}(\tilde{f})_{M, \omega},
\end{equation}
где  ${\cal O}(1)$ --- величины, равномерно ограниченные по $n.$

Нам понадобиться также следующее утверждение, доказанное в работе  \cite{CH_vis_KNU_2013}, которое позволяет находить оценку сверху нормы $(\psi;\beta)$-производной тригонометрического полинома через норму самого полинома в пространстве $L_{M,\omega}$ и  является аналогом  классического результата С.Н.~Бернштейна \cite{Bernshtein-1912}  о  неравенстве для максимума модуля производной тригонометрического полинома. Для его формулирования  через $D^\psi_\beta f$ обозначим $(\psi;\beta)$-производную функции $f,$ a через  $\mathfrak{M}^*$ множество не убывающих и исчезающих на бесконечности последовательностей $\psi(k)>0,$
то есть:
$$
    \mathfrak{M}^*=\{\psi(k):\quad  \psi(k)-\psi(k+1)\ge 0,~
    \lim\limits_{k\to\infty} \psi(k)=0,~ k \in \mathbb{N}\}.
$$

\begin{vsplem}\label{L.B}
Пусть $\psi \in \mathfrak{M}^*,$ $\beta \in \mathbb{R},$ $M \in {\cal Q}_2^\theta$ и $\omega \in A_{p(M)}.$ Тогда для
произвольного тригонометрического полинома $T_n$ порядка, не выше
$n,$ выполняется неравенство
\begin{equation}\label{15}
    \| D_{{\beta}}^\psi T_n\|_{M,\omega}\le \frac{K}{\psi(n)}
    \|T_n\|_{M,\omega}, \quad n=0,1,\ldots,
\end{equation}
где $1/\psi(0):=0,$ a $K$ --- положительная постоянная, которая не зависит от $n$.
\end{vsplem}

Будем также использовать следующее утверждение, доказанное в работе \cite{Akgun-Weighted_Smirnov_Orlicz}.
\begin{vsplem}\label{L.A}
Пусть $M \in {\cal Q}_2^\theta$ и $\omega \in A_{p(M)}.$  Тогда для произвольной функции
$f \in L_{M,\omega}$ имеет место неравенство
\begin{equation}\label{8}
    \Omega_k (f; \delta)_{M, \omega} \le C \delta^{2 k} \| f^{2k}\|_{M,\omega},
    \quad r=1,2,\ldots~.
\end{equation}
\end{vsplem}

Отметим, что исследованию различных вопросов теории приближения в пространствах Орлича посвящено значительное количество работ. Упомянем в этом обзоре
преимущественно результаты, касающиеся получения прямых и обратных теорем теории аппроксимации в этих пространствах.

Так некоторые прямые теоремы в  пространствах Орлича были получены в работах \cite{RAMAZANOV, GARIDI, RUNOVSKI}. Случай весовых пространств Орлича был рассмотрен в работах \cite{Akgun-Weighted_Smirnov_Orlicz, ISRAFILOV_GUVEN}, где был получен ряд прямых и обратных теорем теории приближения в этих пространствах. Некоторые обобщения упомянутых результатов в весовых пространствах Орлича, определенных на различных множествах комплексной плоскости, были доказаны в \cite{KOKILASHVILI_A direct theorem, KOKILASHVILI_On analytic functions, ISRAFILOV_p-Faber polynomials, ISRAFILOV_p-Faber-Laurent}. Кроме того, различные уточнения и обобщения прямых и обратных теорем теории приближения в весовых пространствах Орлича  были доказаны в работах \cite{Akgun_Israfilov_Approx_in_WOS, Guven_Israfilov_On_Approx_in_WOS}.

В этой работе получен ряд  прямых и обратных теорем теории приближения для $\psi$-дифференцируемых и
$(\psi;\beta)$-дифференцируемых функций в  весовых пространствах Орлича $L_{M,\omega}$ с весами, которые принадлежат так называемому классу Макенхаупа. Используемое в этой работе понятие $\psi$-производной обобщает известные понятия $(\psi;\beta)$-производной, дробных производных в смысле Вейля и Вейля-Надя, обычной производной, а полученные результаты, с одной стороны, содержат как частные случаи некоторые результаты предшественников, а с другой стороны, позволяют обнаружить новые факты.


\bigskip

{\bf 3. Основные результаты.}
Пусть $\mathfrak{M}$ ---  подмножество выпуклых последовательностей действительных чисел
$\psi(k)>0,$ из множества $\mathfrak{M}^*,$ то есть:
$$
    \mathfrak{M}:=\{\psi \in \mathfrak{M}^*:~ \psi(k+2)-2\psi(k+1)+\psi(k)>0, \quad k \in \mathbb{N}\},
$$
а $\mathfrak{M}'$ --- подмножество функций $\psi\in \mathfrak{M},$ для
которых  $\sum\limits_{k=1}^\infty \frac{\psi(k)}{k}
<\infty.$

Обозначим через $L^{\psi}L_{M,\omega}$  классы $\psi$-интегралов функций
$f \in  L_{M,\omega}$.  Имеет место следующее утверждение.

\begin{teo}\label{T.1}
Пусть $f \in L^{\psi}L_{M,\omega},$ причем $M \in {\cal Q}_2^\theta,$  $\omega \in A_{p(M)}$
и  $\psi_1 \in \mathfrak{M},~\psi_2 \in \mathfrak{M}'.$ Тогда для произвольного натурального
$n$  выполняется неравенство
\begin{equation}\label{16}
    E_n(f)_{M, \omega}\le K \psi(n)
    E_n(f^\psi)_{M, \omega},
\end{equation}
в котором $\psi(n)=\sqrt{\psi_1^2(n)+\psi_2^2(n)},$ $K$ --- константа, которая не зависит  от $n$ и функции $f.$
\end{teo}

Отметим, что утверждение теоремы \ref{T.1} в большинстве случаев является
хорошо известным. При $\psi(n)=n^{-\alpha},~\beta=\alpha,~
\alpha>0,~n\in \mathbb{N}$ оценка (\ref{16}) получена в работе
\cite{Akgun-Weighted_Smirnov_Orlicz}. При $\psi_1(n)=\psi(n)\sin \frac{\beta\pi}{2},$ $\psi_2(n)=\psi(n)\cos \frac{\beta\pi}{2},$ $\psi \in \mathfrak{M},~\beta \in \mathbb{R}$ соответствующие результаты найдены в работе \cite{CH_vis_KNU_2013}, а если $M= t^p,~1\le p\le\infty,$ то такие
результаты принадлежат А.И.~Степанцу и А.К.~Кушпелю
\cite{Stepanetz-Kushpel-prepr-1984, Stepanetz-Kushpel-umg-1987}.
Если же $1< p<\infty$ и к тому же
$\psi(n)=n^{-r},~\beta=r,~n,r\in \mathbb{N},$ неравенство (\ref{16})
доказано А.Ф.~Тиманом \cite[c.~316]{A.Timan-teoria-priblig}.

\textit{\textbf{Доказательство.}} Пусть $f \in L^{\psi}L_{M,\omega},$ где $M \in {\cal Q}_2^\theta$ и $\omega \in A_{p(M)},$ Тогда, если $\psi_1 \in \mathfrak{M},~\psi_2 \in \mathfrak{M}'$  то \cite[c.~151 -- 152]{Stepanetz-2002m1} (см. также \cite[c.~140]{Stepanetz-2002-angl}), всегда
существует функция $f^\psi \in L^0,~L^0:=\{\varphi \in L:~\varphi \perp 1\},$ ряд Фурье
которой совпадает с (\ref{3}). В работе \cite[лемма 3]{Khabazi} показано, что для данной функции $f \in L_{M, \omega}$ и произвольного $\varepsilon > 0$ всегда найдется тригонометрический полином $T(x),$ для которого
$$
    \int\limits_{0}^{2\pi} M( |f(x)-T(x)| ) \omega(x)~dx<\varepsilon.
$$
Отсюда, в частности, следует, что
$$
    E_n(f)_{M,\omega} \to 0, \quad n \to \infty.
$$

Поэтому, учитывая соотношение (\ref{7}), можем утверждать, что в смысле сходимости в метрике пространств $L_{M,\omega}$ выполняются равенства
\begin{equation}\label{11}
        f(x)=\sum_{k=0}^\infty A_k (f; x), \quad
    f^\psi(x)=\sum_{k=1}^\infty A_k (f^\psi; x)=
$$
$$
    =\sum_{k=1}^{\infty}\frac{\psi_1(k)}{\psi^2(k)}
    A_k(f;x)-\frac{\psi_2(k)}{\psi^2(k)} \tilde{A}(f;x).
\end{equation}

Из разложения (\ref{11}) следуют формулы связи между коэффициентами
Фурье функций $f^\psi$ и $f$:
\begin{equation}\label{14.1}
    a_k(f)=\psi_1(k) a_k(f^\psi)- \psi_2 (k) b_k (f^\psi),
\end{equation}
\begin{equation}\label{15.1}
    b_k(f)=\psi_2(k) a_k(f^\psi)+ \psi_1(k) b_k (f^\psi).
\end{equation}

Учитывая теперь равенства (\ref{14.1}) -- (\ref{15.1}), находим
$$
    f(x)=\sum_{k=0}^\infty A_k (f; x) =\frac{a_0(f)}{2}+
    \sum_{k=1}^\infty(a_k(f)\cos~kx+b_k(f)\sin~kx)=
$$
$$
    =\frac{a_0(f)}{2}+\sum_{k=1}^\infty \bigg( \bigg[\psi_1(k) a_k(f^\psi)-
    \psi_2 (k) b_k (f^\psi) \bigg]\cos kx+
$$
$$
    + \bigg [\psi_2(k) a_k(f^\psi)+
    \psi_1(k) b_k (f^\psi) \bigg] \sin kx \bigg)=
$$
$$
    =\frac{a_0(f)}{2}+\sum_{k=1}^\infty \bigg(\psi_1 (k)
    \bigg [a_k(f^\psi) \cos kx+
    b_k (f^\psi) \sin kx \bigg] +
$$
$$
    +\psi_2(k) \bigg [a_k(f^\psi)\sin kx-
    b_k (f^\psi) \cos kx \bigg]  \bigg)=
$$
\begin{equation}\label{17}
    =\frac{a_0(f)}{2} + \sum_{k=1}^\infty \psi_1 (k) A_k(f^\psi;x)
    +\sum_{k=1}^\infty \psi_2 (k) \tilde{A}_k(f^\psi;x).
\end{equation}

Принимая во внимание соотношения (\ref{9}) и  (\ref{17}), получаем
$$
    f(x)-S_{n-1}(f;x)=\sum_{k=n}^\infty A_k(f;x)=
$$
\begin{equation}\label{18}
    =\sum_{k=n}^\infty \psi_1(k) A_k(f^\psi;x)  +
    \sum_{k=n}^\infty \psi_2 (k) \tilde{A}_k(f^\psi;x).
\end{equation}

Далее имеем
$$
    \sum_{k=n}^\infty \psi_1  (k) A_k(f^\psi;x)=
    \sum_{k=n}^\infty \psi_1  (k) \Bigg[ \bigg( S_{k}(f^\psi; x)- f^\psi(x)\bigg)-
    \bigg( S_{k-1}(f^\psi;x)- f^\psi(x)\bigg)\Bigg]=
$$
$$
    =\sum_{k=n}^\infty \psi_1  (k) \bigg(  S_{k}(f^\psi; x)- f^\psi(x)\bigg)-
    \sum_{k=n}^\infty \psi_1  (k) \bigg( S_{k-1}(f^\psi;x)-f^\psi(x)\bigg)=
$$
\begin{equation}\label{19}
    =\sum_{k=n}^\infty [\psi_1  (k)-\psi_1  (k+1) ]
    \bigg[S_{k}(f^\psi; x)-f^\psi(x)\bigg]-\psi_1 (n)
    \bigg[S_{n-1}(f^\psi; x)-f^\psi(x)\bigg].
\end{equation}

Аналогично получаем
$$
   \sum_{k=n}^\infty \psi_2  (k) \tilde{A}_k(f^\psi;x)= \sum_{k=n}^\infty \psi_2  (k) A_k(\tilde{f}^\psi;x)=
$$
\begin{equation}\label{20}
   =\sum_{k=n}^\infty [\psi_2  (k)-\psi_2  (k+1) ]
    \bigg[S_{k}(\tilde{f}^\psi; x)-
    \tilde{f}^\psi(x)\bigg]-
    \psi_2 (n) \bigg[S_{n-1}(\tilde{f}^\psi; x)-\tilde{f}^\psi(x)\bigg].
\end{equation}

Используя свойства нормы и соотношение (\ref{9}), на
основании равенств (\ref{18}) -- (\ref{20}) находим
$$
    E_n (f)_{M, \omega}\le \| f-S_{n-1}(f) \|_{M,\omega} \le
$$
$$
    \le \sum_{k=n}^\infty [\psi_1  (k)-\psi_1  (k+1) ]
    \bigg\| S_{k}(f^\psi)-f^\psi \bigg\|_{M,\omega}+\psi_1 (n)
    \bigg\| S_{n-1}(f^\psi)- f^\psi \bigg\|_{M,\omega}+
$$
$$
    +\sum_{k=n}^\infty [\psi_2  (k)-\psi_2  (k+1) ]
    \bigg\| S_{k}(\tilde{f}^\psi)-
    \tilde{f}^\psi \bigg\|_{M,\omega}+
    \psi_2 (n) \bigg\|S_{n-1}(\tilde{f}^\psi)-
    \tilde{f}^\psi \bigg\|_{M,\omega} \le
$$
$$
    \le K \bigg(\psi_1  (n) E_n (f^\psi)_{M, \omega}+
    \psi_2  (n) E_n(\tilde{f}^\psi)_{M, \omega}\bigg)\le
    K \psi (n) E_n (f^\psi)_{M, \omega},
$$
где  $\psi^2(n)=\psi_1^2(n)+\psi_2^2(n)$ и $K$ --- константа, которая не зависит  от  $n.$  Теорема доказана.

\begin{sled}\label{s1}
В условиях теоремы 1 имеет место  неравенство
\begin{equation}\label{16*}
    E_n(f)_{M, \omega}\le K \psi(n),
\end{equation}
в котором $K$ --- константа, которая не зависит от $n.$
\end{sled}

Считая, что последовательности $\psi(k)$ из  $\mathfrak{M}$ являются
сужениями на множество натуральних чисел непрерывных функций
$\psi(t)$ непрерывного аргумента $t\ge 1,$ в соответствии с
\cite[c. 159]{Stepanetz-2002m1} (см. также \cite[c.~147]{Stepanetz-2002-angl}) через $\eta(t)=\eta(\psi;t)$
обозначим функцию, которая связана с $\psi$ равенством
$$
    \psi(\eta(t))=\frac{1}{2} \psi(t), \quad t\ge 1.
$$
Отсюда, вследствие строгой монотонности и убывания к нулю $\psi,$
функция $\eta(t)$ для всех $t\ge 1$ определяется однозначно
\begin{equation}\label{26.1}
    \eta(t)=\eta(\psi;t)=\psi^{-1}\bigg(\frac{1}{2} \psi(t)\bigg),
\end{equation}
где $\psi^{-1}$ --- функция, обратная к $\psi.$

Для получения следующих результатов будем  использовать определения множеств
$\mathfrak{M}_0$ и $F,$ принадлежащие А.И.~Степанцу \cite[c. 160, c. 165]{Stepanetz-2002m1} (см. также \cite[c.~148, 153]{Stepanetz-2002-angl}):
$$
\mathfrak{M}_0=\{\psi \in \mathfrak{M}: \quad
0<\frac{t}{\eta(\psi;t)-t}\le K, ~t\ge1 \},
$$
$$
F=\{\psi \in \mathfrak{M}: \quad \eta'(\psi; t) \le K   \},
$$
где $\eta(\psi;t)$ --- функция, которая определяется равенством (\ref{26.1}).

Справедлива следующая теорема.

\begin{teo}\label{T.2}
Пусть $M \in {\cal Q}_2^\theta$ и $\omega \in A_{p(M)}.$  Тогда для произвольной функции $f \in L_{M,\omega}$ и $n\in \mathbb{N}$ будем иметь:

\begin{enumerate}
\item
Если $\psi\in \mathfrak{M}_0$  и ряд
$$
\sum_{k=1}^\infty E_k(f)_{M, \omega} |k\psi(k)|^{-1},
$$
сходится,  тогда существует производная $f^\psi_\beta,$ такая что
для $r=0,1,2,\dots,$ и произвольного натурального $n$  выполняется неравенство
\begin{equation}\label{21}
    \Omega_r \bigg(f^\psi_\beta; \frac{1}{n} \bigg)_{M, \omega} \le
    \frac{C}{n^{2r}} \sum\limits_{\nu=0}^{n} \frac{\nu^{2r}}{ \psi(\nu)}E_{\nu}(f)_{M,\omega}+
    C\sum_{\nu=n+1}^\infty \frac{E_{\nu}(f)_{M, \omega} }{\nu \psi(\nu)},
\end{equation}
\item
Если $\psi\in F,~ \eta(\psi;t)-t\ge C>0$ и ряд
$$
    \sum_{k=1}^\infty E_k(f)_{M, \omega} |\psi(k)(\eta(k)-k)|^{-1},
$$
сходится, то существует производная $f^\psi_\beta$ для которой
$$
    \Omega_r \bigg(f^\psi_\beta; \frac{1}{n}\bigg)_{M, \omega} \le
    C \Bigg( \frac{1}{n^{2r}}
    \sum_{\nu=1}^{n} \frac{\nu^{2r}}{\psi(\nu) } E_\nu(f)_{M, \omega}+
$$
\begin{equation}\label{41.3}
    + \sum_{\nu=n+1}^{\infty} \frac{ E_\nu (f)_{M, \omega}}{\psi(\nu)(\eta(\nu)-\nu)}\Bigg),
    \quad n \in \mathbb{N},
\end{equation}

\end{enumerate}
где  $C$ --- константа, которая не зависит  от $n$  и функции $f.$
\end{teo}

Если $\psi(n)=n^{-r},~\beta=r,~ r,~n\in
\mathbb{N},$ тогда утверждение теоремы  \ref{T.2} совпадает с теоремой  1 из работы \cite{Guven_Israfilov_On_Approx_in_WOS}. В случае  $\psi(n)=n^{-\alpha},~\beta=\alpha,~ \alpha>0,~n\in
\mathbb{N},$ аналогичные результаты были получены в работе \cite{Akgun-Weighted_Smirnov_Orlicz}.

\textit{\textbf{Доказательство.}} Для доказательства теоремы будем использовать схему, предложенную в книге \cite[c. 120 -- 126]{Stepanetz-2002m2} (см. также \cite[c.~545 -- 551]{Stepanetz-2002-angl}). Пусть $\{t_n(\cdot) \}_{n=1}^\infty$
последовательность тригонометрических полиномов, которые осуществляют наилучшее приближение функции $f\in L_{M,\omega}.$
Положим для данной функции $\psi$ и каждого натурального $n$
$$
    n_0=n, n_1=[\eta(\psi; n)]+1,\ldots,n_i=[\eta(\psi; n_{i-1})]+1, \ldots~.
$$
В этом случае ряд
$$
    t_{n_0}(x)+\sum_{i=1}^\infty (t_{n_i}(x)-t_{n_{i-1}}(x))
$$
будет сходиться к функции $f$ в пространстве $L_{M, \omega}$ . Рассмотрим ряд
\begin{equation}\label{22}
    (D^\psi_\beta t_{n_0}) (x) +\sum_{i=1}^\infty
    (D^\psi_\beta[t_{n_i}-t_{n_{i-1}}])(x)
\end{equation}
и убедимся, что он будет сходиться  в пространстве $L_{M, \omega}$ к сумме $T(x),$
ряд Фурье которой имеет вид
\begin{equation}\label{23}
    \sum_{k=1}^{\infty}\frac{1}{\psi(k)}
    \left( a_k(f)\cos\bigg(kx+\frac{\beta \pi}{2}\bigg)+
    b_k(f)\sin\bigg(kx+\frac{\beta \pi}{2}\bigg) \right).
\end{equation}

Применяя неравенство (\ref{15}) к разности  $u_i(x)=t_{n_i}(x)-t_{n_{i-1}}(x),$
(которая, очевидно, является полиномом порядка  $n_i$), получаем
$$
    \| (D^\psi_\beta u_i)(\cdot) \|_{M,\omega} \le
    C E_{n_{i-1}+1}(f)_{M, \omega} |\psi(n_i)|^{-1},
$$
вследствие чего
$$
    \sum_{i=1}^\infty \| (D^\psi_\beta u_i)(\cdot) \|_{M,\omega}
    \le C\bigg(E_{n+1}(f)_{M, \omega} (\psi(n))^{-1}+\\
$$
\begin{equation}\label{24}
    +\sum_{i=1}^\infty E_{n_i+1}(f)_{M, \omega}|\psi(n_i)|^{-1} \bigg).
\end{equation}

Используя оценку
\begin{equation}\label{31}
    \frac{E_{n_i+1}(f)}{\psi(n_i)}\le \sum_{\nu=n_{i-1}}^{n_i-1}
    \frac{E_{\nu+1}(f)}{\nu \psi(\nu)}, \quad \psi \in \mathfrak{M}_0,
\end{equation}
полученную в книге \cite[c. 124 -- 125]{Stepanetz-2002m2} (см. также \cite[c.~149]{Stepanetz-2002-angl}), из соотношения (\ref{24}) получаем
\begin{equation}\label{25}
    \sum_{i=1}^\infty \| (D^\psi_\beta u_i)(\cdot) \|_{M,\omega}
    \le C\bigg(\frac{E_{n+1}(f)_{M, \omega}}{ \psi(n)}+
    \sum_{\nu=n+1}^\infty \frac{E_{\nu}(f)_{M, \omega}}{\nu \psi(\nu)} \bigg).
\end{equation}

Поскольку по условию  теоремы ряд в правой части соотношения (\ref{25}) сходится, то
это означает, что ряд (\ref{23}) действительно сходится в пространстве $ L_{M,\omega}$ к некоторой функции $T(x)\in  L_{M,\omega}.$

Пусть $a_k^{(n_i)}=a_k (t_{n_i})$ и $b_k^{(n_i)}=b_k (t_{n_i}),~
k=0,1,2,\ldots,$ --- коэффициенты Фурье полиномов $t_{n_i}(x).$
Тогда в соответствии с равенствами  (\ref{14.1}) и (\ref{15.1}) коэффициенты
$\alpha_k^{(n_i)}$ и $\beta_k^{(n_i)}$ полиномов $(D^\psi_\beta t_{n_i})(\cdot)$ имеют вид
\begin{equation}\label{26}
    \alpha_k^{(n_i)}=\frac{1}{\psi(k)} \bigg (a_k^{(n_i)}\cos \frac{\beta  \pi}{2}+
    b_k^{(n_i)}\sin \frac{\beta  \pi}{2} \bigg)
\end{equation}
\begin{equation}\label{27}
    \beta_k^{(n_i)}=\frac{1}{\psi(k)} \bigg (b_k^{(n_i)} \cos \frac{\beta  \pi}{2}
    - a_k^{(n_i)} \sin \frac{\beta  \pi}{2}\bigg).
\end{equation}

Поскольку равенство
$$
 T(x)=\lim_{n\to \infty} (D^\psi_\beta t_{n_i})(x)
$$
выполняется в смысле сходимости в метрике пространств $ L_{M,\omega},$ то
$$
    a_k(T)=\lim_{i\to \infty} \alpha_k^{(n_i)}, \quad
    b_k(T)=\lim_{i\to \infty} \beta_k^{(n_i)}, \quad k=0,1,\ldots~.
$$
Принимая во внимание то, что
$$
    \lim_{i\to \infty} a_k^{(n_i)}=a_k(f), \quad
    \lim_{i\to \infty} b_k^{(n_i)}=b_k(f), \quad k=0,1,\ldots~,
$$
из равенств (\ref{26}) -- (\ref{27}) получаем
$$
    a_k(T)=\frac{1}{\psi(k)} \bigg (a_k(f) \cos \frac{\beta  \pi}{2}+
    b_k(f) \sin \frac{\beta  \pi}{2} \bigg)
$$
$$
    b_k(T)=\frac{1}{\psi(k)} \bigg (b_k(f) \cos \frac{\beta  \pi}{2}
    - a_k(f) \sin \frac{\beta  \pi}{2}\bigg).
$$
Отсюда следует, что ряд Фурье функции $T(x)$ действительно совпадает с рядом (\ref{23}). Это означает, что функция $f(x)$ действительно имеет $(\psi;\beta)$-производную  $f^\psi_\beta(x),$ которая принадлежит пространству $L_{M, \omega}$, и для неё выполняется равенство
\begin{equation}\label{28}
    f^\psi_\beta(x)= (D^\psi_\beta t_{n_0}) (x) +
    \sum_{i=1}^\infty (D^\psi_\beta[t_{n_i}-t_{n_{i-1}}])(x),
\end{equation}
в метрике пространства $ L_{M, \omega}.$

Используя свойство 2 величины $\Omega_r (f^\psi ; \frac{1}{n})_{M, \omega}$, на основании соотношения (\ref{28}) получаем
\begin{equation}\label{29}
    \Omega_r \bigg(f^\psi ; \frac{1}{n}\bigg)_{M, \omega} \le
    \Omega_r \bigg(f^\psi - D^\psi_\beta t_{n_0}]; \frac 1n\bigg)
    +\Omega_r \bigg(D^\psi_\beta t_{n_0}; \frac{1}{n}\bigg)_{M, \omega}
\end{equation}

Учитывая свойство 3 модуля непрерывности и неравенство (\ref{15}), находим
$$
    \Omega_r \bigg(f^\psi - D^\psi_\beta t_{n_0}; \frac{1}{n}\bigg)_{M,\omega} \le
    C \|f^\psi - D^\psi_\beta t_{n_0}] \|_{M,\omega}=
$$
\begin{equation}\label{30}
    =C\| \sum\limits_{i=1}^\infty (D_\beta^\psi (t_{n_i}- t_{n_{i-1}}))(x)\|_{M,\omega}
    \le  C \sum\limits_{i=1}^\infty \frac{1}{\psi(n_i)} \| t_{n_i}- t_{n_{i-1}}\|_{M,\omega}\le
$$
$$
    \le  C\sum\limits_{i=1}^\infty \frac{\| t_{n_i}-f \|_{M,\omega}}{\psi(n_i)}  + C\sum\limits_{i=1}^\infty \frac{\|f-t_{n_{i-1}}\|_{M,\omega}}{\psi(n_i)} \le
    K \sum\limits_{i=1}^\infty \frac{E_{n_i+1}(f)_{M, \omega}}{\psi(n_i)}.
\end{equation}

Используя оценку (\ref{31}) из соотношения (\ref{30}) получаем
\begin{equation}\label{32}
    \Omega_r \bigg(f^\psi - D^\psi_\beta t_{n_0}; \frac{1}{n}\bigg)_{M,\omega} \le
    C \sum_{\nu=n_0+1}^\infty \frac{E_{\nu}(f)_{M, \omega} }{\nu \psi(\nu)}.
\end{equation}

Принимая во внимание неравенства (\ref{8}) и (\ref{15}), находим
$$
    \Omega_r \bigg(D^\psi_\beta t_{n_0}; \frac{1}{n} \bigg)_{M, \omega} \le
    \Omega_r \bigg( D^\psi_\beta t_{0}; \frac{1}{n} \bigg)_{M, \omega}+
    \sum\limits_{i=1}^{n_0} \Omega_r \bigg (D_\beta^\psi (t_{i}-
     t_{i-1});\frac 1n)_{M,\omega}\le
$$
$$
    \le \frac{C}{n^{2r}} \bigg(  \bigg\|(D^\psi_\beta t_{0})^{(2r)} \bigg\|_{M, \omega}+
    \sum\limits_{i=1}^{n_0} \bigg \| (D_\beta^\psi (t_{i}-
     t_{i-1}))^{(2r)} \bigg\|_{M,\omega}\bigg)\le
$$
\begin{equation}\label{33}
    \le \frac{C}{n^{2r}}
    \sum\limits_{\nu =0}^{n_0}  \nu^{2r} \frac{E_{\nu}(f)_{M,\omega}}{\psi(\nu)}.
\end{equation}

Объединяя теперь соотношения (\ref{29}), (\ref{32}) и (\ref{33}), получаем
$$
    \Omega_r \bigg(f^\psi_\beta; \frac{1}{n} \bigg)_{M, \omega} \le
    \frac{C}{n^{2r}} \sum\limits_{\nu=0}^{n_0} \frac{\nu^{2r}}{ \psi(\nu)}E_{\nu}(f)_{M,\omega}+
    C\sum_{\nu=n_0+1}^\infty \frac{E_{\nu}(f)_{M, \omega} }{\nu \psi(\nu)},
$$

и, поскольку по определению   $n_0 = n,$ заканчиваем доказательство теоремы.

Доказательство пункта 2 теоремы проводится аналогично, с учетом следующего аналога неравенства (\ref{31}), полученного в книге \cite[c.~125 -- 126]{Stepanetz-2002m2} (см. также \cite[c.~551]{Stepanetz-2002-angl})
$$
    \frac{E_{n_i+1}(f)}{\psi(n_i)}\le \sum_{\nu=n_{i-1}}^{\eta(n_{i-1})}
    \frac{E_{\nu+1}(f)}{(\eta(\nu)-\nu) \psi(\nu)}, \quad \psi \in F.
$$



\bibliographystyle{plain}
\renewcommand{\refname}{Литература}

\end{document}